\documentclass[12pt]{article}

\usepackage{graphicx}
\usepackage{caption}

\begin{document}

\begin{center}
{\bf\large To the Memory of Alexander Grothendieck:}

\vspace*{3mm}

{\bf\large a Great and Mysterious Genius of Mathematics} \\

\vspace*{5mm}

{\small Wolfgang Bietenholz$^{\rm a}$ and Tatiana Peixoto$^{\rm b}$ \\

\vspace*{2mm}

$^{\rm a}$ Instituto de Ciencias Nucleares \\
Universidad Nacional Aut\'{o}noma de M\'{e}xico \\
A.P. 70-543, C.P. 04510 Distrito Federal , Mexico

\vspace*{2mm}

$^{\rm b}$ Universidade Federal do ABC \\
C\^{a}mpus Santo Andr\'{e}, CEP 09606-070, Brazil

}

\end{center}

{\small In November 2014 Alexander Grothendieck passed away at the age 
of 86. There is no doubt that he was one of the greatest and
most innovative mathematicians of the 20th century. 
After a bitter childhood, his meteoric ascent 
started in the Cartan Seminar in Paris, it led to a breakthrough while
he worked in S\~{a}o Paulo, and to the Fields Medal.
He introduced numerous new concepts and techniques,
which were involved in the groundbreaking solutions to long-standing
problems. However, dramatic changes were still ahead of him.
In recent years hardly anybody knew where he was living, and even 
if he was still alive; he had withdrawn to a modest life in isolation.
Also beyond his achievements in mathematics, Grothendieck was an
extraordinary person. This is a tribute of his fascinating life.}

\vspace*{2mm}

\begin{center}
{\bf\large Em mem\'{o}ria de Alexander Grothendieck:}

\vspace*{3mm}

{\bf\large um grande e misterioso g\^{e}nio da matem\'{a}tica} \\
\end{center}

{\small Em Novembro de 2014 Alexander Grothendieck faleceu aos seus 86 anos 
de idade. N\~{a}o h\'{a} duvidas de que ele foi um dos maiores e mais 
criativos matem\'{a}ticos do s\'{e}culo XX. Ap\'{o}s uma inf\^{a}ncia amarga, 
sua ascens\~{a}o mete\'{o}rica iniciou no Semin\'{a}rio Cartan em Paris, o 
que o levou a um avan\c{c}o enquanto ele trabalhava em S\~{a}o Paulo, 
e a Medalha Fields. Ele introduziu numerosos novos conceitos e 
t\'{e}cnicas, que foram envolvidos nas solu\c{c}\~{o}es 
inovadoras de problemas de longa data. No entanto, 
mudan\c{c}as dram\'{a}ticas ainda estavam por vir. 
Nos \'{u}ltimos anos, quase ningu\'{e}m sabia onde ele estava morando, 
ou at\'{e} mesmo se ele estava vivo; ele havia se retirado para uma 
vida modesta em isolamento. Tamb\'{e}m para al\'{e}m dos seus m\'{e}ritos 
em matem\'{a}tica, Grothendieck foi uma pessoa extraordin\'{a}ria. 
Isso \'{e} um tributo a sua fascinante vida.}

\newpage

\section{An unusual family}

Alexander Grothendieck's life was dominated by turbulence and
radical turning points. As constant features, however, he followed 
consistently his own path --- keeping away from anything established ---
and whatever he did, he did with absolute passion.

To capture the spirit of his highly unusual biography,\footnote{We adopt
the biographical information essentially from Refs.\ \cite{Jackson,Scharlau}, 
but the most complete source is Ref.\ \cite{Scharlau2}. We apologize
for not quoting them each time.}
we have to start with his parents.
His father --- whom Grothendieck honored very much --- was 
Alexander Schapiro (1890-1942), born in the Russian town
Novozybkov (in the border region with Belarus and Ukraine), 
in an orthodox Jewish community. Still very young he joined an 
armed anarchist group, which was captured in 1905, after the 
failed attempt to overthrow the tsarist regime. All members
were executed, except for Alexander, who was pardoned to
life in prison due to his youth. About ten years later he escaped,
and readily joined another anarchist army, this time in the Ukraine.
He was captured again, sentenced to death, but he managed
to escape once more (though he lost his left arm).

Then he lived under the name Alexander Taranow in Berlin 
(and other cities), where he worked as an independent photographer. 
Around 1924 he met
Hanka (actually Johanna) Grothendieck (1900-57), told her husband:
``I will steal your wife'' --- and he proceeded in doing so.
Hanka was a far-left activist too, and she tried to become a 
journalist and writer, 
but --- despite her  talent --- she could not publish much.
In 1928 their son Alexander
Grothendieck was born; Schurik --- this was his nickname --- lived 
for the first five years with his parents and a half-sister in 
Berlin (in the Scheunenviertel).

\section{Youth during World War II}

In 1933, when the Nazis came to power, the situation was getting
too dangerous for Schurik's father, who flew to Paris. Hanka joined
him soon, and left Schurik with a foster family in Hamburg.
There he attended school from 1934-39 and lived in the home 
of Wilhelm Heydorn --- a former military
officer and priest, who turned pacifist and atheist. 
Alexander's Jewish ancestry was kept secret, but in 1939 Germany 
was getting too dangerous for him as well, in particular because
his foster parents opposed the Nazi regime. He was put on a train
to France where he met his parents again; they were back from the 
Spanish civil war, having supported an anarchist group.
In 1940 the family was imprisoned in internment camps 
by the Vichy regime, which collaborated with the Germans.
Two years later, Alexander Taranow was extradited to the Nazis, 
he was deported to Auschwitz where he died.

\begin{figure}
\begin{center}
\includegraphics[angle=0,width=.27\linewidth]{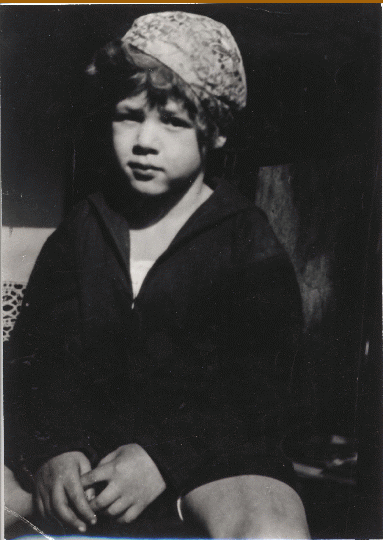}
\end{center}
\vspace*{-2mm}
\caption*{\it Alexander (Schurik) Grothendieck a the age of 12.} 
\vspace*{-2mm}
\end{figure}

Also in 1942, 14-year-old Alexander arrived in
Le Chambon-sur-Lignon, a small town in the Massif Central, which
was a center of resistance against the German occupation.
Here Alexander attended a school, which was devoted to the
spirit of pacifism. When there were raids by the Gestapo, he and 
other pupils hided in the forest for a couple of days,
divided into small groups \cite{Recoltes}. 
In 1945 he finished his {\it baccalaureate.}

After the war his mother Hanka was released, and Alexander was
closely attached to her until her death in 1957.\footnote{Hanka
died of tuberculosis, probably as a consequence of her confinement
during the war. 
However, she still witnessed her son's upcoming world fame.} 
They moved to Montpellier, where Alexander 
studied mathematics, and received a modest scholarship.
The local university was not very helpful to him,
so he resorted mostly to auto-didactic studies. He was particularly 
interested in a deep understanding of space and geometry, starting 
with the notion of a {\em point,}\footnote{Already in High School it upset 
him that his text books never gave a satisfactory definition of terms 
like {\em length, area} and {\em volume.}} and he elaborated by himself 
a generalized concept of integration.\footnote{Apparently 
Hanka obstructed his contact with girls, which was certainly 
favorable for the intensity of his studies. Nevertheless, in his 
later life he had five children with three women. One of them, Mireille 
Dufour, was his wife in the 1960s. She was from the Normandy,
a little older than Alexander, and she also had links to the anarchist
movement in Spain. It is reported that even their marriage had an
touch of anarchy; now Alexander had a variety of affairs.\label{women}}

\section{A fairytale-like career}

In 1948 Alexander was awarded a fellowship to go to Paris, where he 
got in contact with mathematical research; in particular, 
he attended the famous {\em Cartan Seminar.} He was not shy
to discuss with famous scholars, he was ambitious and 
passionate; later he wrote ``j'\'{e}tais un 
math\'{e}maticien: quelqu'un qui {\em fait} des maths, au plein 
sens du terme -- comme on {\em fait} l'amour'' \cite{Recoltes}.
Initially he hoped for his independent work to
provide a quick Ph.D., but he was told that he had
essentially re-discovered the Lebesgue integral (which
had been known since the early 20th century).
Also later, as a highly established mathematician, he
always followed his own ideas, rather than studying the 
literature (he got informed about relevant results
in discussions).

Since Alexander wanted to explore Topological Vector Spaces, 
Henri Cartan and Andr\'{e} Weil recommended him to move to the
University of Nancy, in Northern France, where two leading
experts were working: Jean Dieudonn\'{e} and Laurent Schwartz
--- the latter was also a pioneer in Distribution Theory,
and he just won a Fields Medal.{\footnote{The Fields 
Medal is the highest distinction in mathematics, sometimes regarded 
equivalent to a Nobel Prize (though there is an age limit of 40). 
It is awarded at the {\it International Congress of Mathematicians} 
(ICM), the largest mathematical conference, which is held once every 
four years (it coincides with the year of the Soccer World Cup).}
He showed his new student his latest paper; it ended with a
list of 14 open questions, relevant for locally convex spaces.
Alexander went ahead and introduced new methods, which allowed
him to solve {\em all} these problems within a few months! A 
mathematical superstar appeared, at the age of 22, with a 
chaotic youth and sparse education.

Despite his success, it was difficult for him to find a job in
France, in particular because he was stateless.
His advisors mentioned the unfortunate situation of this 
young genius, and found a visitor position for him at the
University of S\~{a}o Paulo,
where Alexander stayed from 1952-54 \cite{Brasil}. The contact was 
established by Paulo Ribenboim, a Brazilian student
of the same age as Alexander, who also worked in Nancy (later
he became a prominent mathematician in Canada). In that period, 
legendary president Get\'{u}lio Vargas was in power in Brazil,
and Alexander finished his Ph.D.\ thesis on {\it Tensor Products 
and Nuclear Spaces} (the second term he had introduced himself).
According to Dieudonn\'{e}, at that time he had already results 
that would have been sufficient for six theses, covering also
functional analysis. He published in Brazilian journals (in French), 
where he introduced the Grothendieck Inequality, 
and he lectured on Topological Vector Spaces; his lecture
notes were published as well \cite{lectTVS}. 
Meanwhile he started to shift his focus of interest towards 
{\em Algebraic Geometry} --- the field where he ultimately had his 
strongest impact; it involves the systematic analysis of the geometric
properties of the solutions to polynomial equations.
\begin{figure}
\begin{center}
\includegraphics[angle=0,width=.35\linewidth]{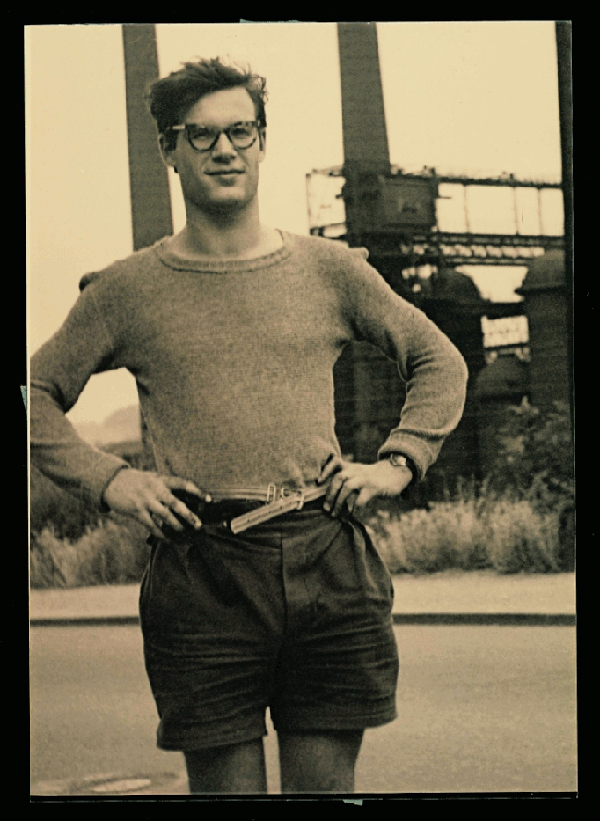}
\end{center}
\vspace*{-2mm}
\caption*{\it The young mathematician Grothendieck, in a picture 
taken by Paulo Ribenboim in 1951.} 
\vspace*{-2mm}
\end{figure}

Alexander worked with full intensity, 
essentially he only paused to sleep and eat.
Schwartz, his Ph.D.\ advisor, asked Ribenboim to encouraged him 
to do occasionally other things in life --- at that time without success.
His mother visited him in Brazil, which presumably again prevented
him from being side-tracked (cf.\ footnote \ref{women}).
His colleague Chaim H\"{o}nig 
(who came to Brazil before the war, as a German refugee)
later remembered that Alexander led a 
``spartan and lonely existence'', living sometimes of milk and 
bananas, and he got frustrated when he failed to solve a problem 
despite hard work.
Still, the problems that he did solve, and the methods that he 
introduced, boosted his meteoric career.

\section{The Golden Era at IH\'{E}S, 1958-70}

After a short stay in Kansas, Grothendieck 
returned to France. Together with Dieudonn\'{e} and
Jean-Pierre Serre, he soon worked at the newly founded 
{\it Institut des Hautes \'{E}tudes Scientifiques} (IH\'{E}S) near 
Paris, which became famous for its research in mathematics
and theoretical physics. Grothendieck
led a group of brilliant young mathematicians.
This era of excellence, 1958-70, coincided with the climax of the 
{\em Bourbaki} group, which Grothendieck was in contact with
(for some years he was a member).
When a visitor noticed that the library of the new institute
was rather incomplete, Grothendieck replied: ``We don't read
books, we write them'' \cite{Jackson}.

His former colleague Pierre Cartier asserts that he run ``one
of the most prestigious mathematics seminars that the world
has ever seen'' \cite{Cartier}. It attracted top mathematicians
from France and all over the world. 
Session could take 10 to 12 hours, leading to improvised
notes that Grothendieck gave to Dieudonn\'{e}, who would then
rewrite them in a neat form. Grothendieck is remembered as an
excellent teacher, who explained also ``trivial'' points patiently, 
with a talent to suggest the appropriate subject to 
each member of his group.
His motivation was simply to {\em understand,} not competition. 

He had anticipated his research program for these years in a
plenary talk at the {\it International Congress of Mathematicians}
in Edinburgh, 1958. His style was to search for ever 
increasing {\em generality and abstraction} (which was a
trend of mathematics in the 20th century), introducing 
accurate new terms and concepts, and working out their properties.  
His colleague John Tate emphasizes that Grothendieck found again and 
again exactly the right level of abstraction, so he was neither dealing
with a special case, nor with a pointless ``vacuum''.
This led to thousands of pages on the merger of Algebraic
Geometry, Arithmetics and Topology. His interest was mostly in 
new, generic concepts, like {\it schemes, \'{e}tales, toposes} 
and {\it motives}; for popular descriptions
we refer to Refs.\ \cite{Cartier,MumTat}, or (more detailed and
technical) Ref.\ \cite{RisingSea}. Grothendieck hardly appreciated 
applications in natural science, like physics.\footnote{Occasionally he got 
a bit interested in biology, encouraged by a friend in Romania.} Even 
the proofs of explicit mathematical theorems were an inspiration for him, 
but not really the ultimate goal. However, Gerd Faltings' proofs of the
Tate and the Mordell Conjecture, as well as Andrew Wiles' proof of 
Fermat's Last Theorem, can all be viewed as applications of {\it motives.}

For a number of years, Grothendieck proved step by step aspects
of the Weil Conjectures (dating back to 1949), which inspired 
amazing new concepts. Later, in 1974, his former student Pierre Deligne 
proved the last point of these conjectures. However, he invoked a 
classical result, deviating from the program of a generalized context,
which employs {\it motives,} as sketched in the IH\'{E}S seminars. His 
mentor appreciated this success, but he was still somehow disappointed.
\begin{figure}
\begin{center}
\includegraphics[angle=0,width=.5\linewidth]{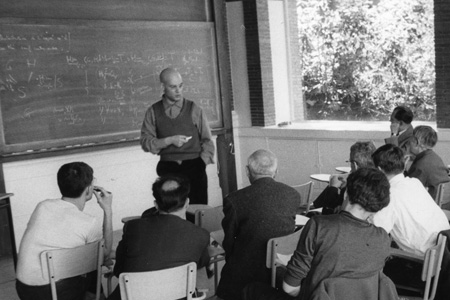}
\hspace*{1cm}
\includegraphics[angle=0,width=.27\linewidth]{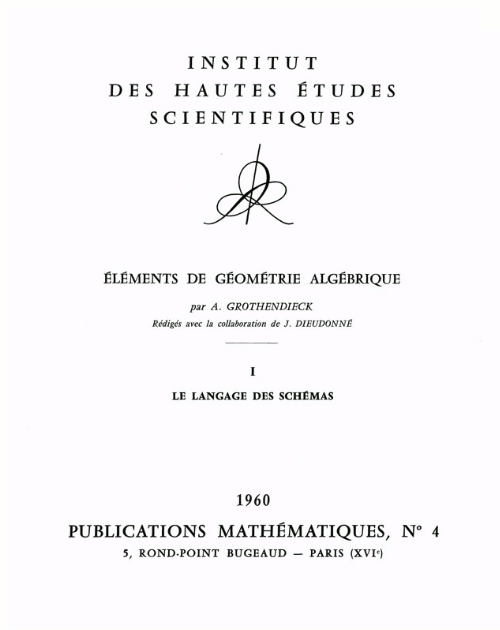}
\end{center}
\vspace*{-2mm}
\caption*{\it On the left: A seminar at the {\it Institut des 
Hautes \'{E}tudes Scientifiques} (IH\'{E}S) near Paris.
On the right: front page of Grothendieck's book
{\it ``\'{E}l\'{e}ments de G\'{e}om\'{e}try Alg\'{e}brique''.}}
\vspace*{-2mm}
\end{figure}

In 1970, at the age of 42, Grothendieck abruptly resigned from 
IH\'{E}S, and entered a completely different phase of his life. 
Soon afterwards his brilliant research group fell apart.

\section{New interests, and a new life style}

Until 1970 Grothendieck's life was almost non-stop focused on 
mathematics (later he called it his ``long period of mathematical
frenzy'' \cite{Recoltes}), and his life style seemed rather bourgeois.
People described him as friendly, direct, by no means arrogant,
idealistic but --- for issues beyond mathematics --- somewhat
na\"{\i}ve. However, other issues of the world did come to his mind,
and gradually became dominant. In particular, he felt strongly
committed to {\em pacifism.}

Since the late 1950s he was wearing Russian peasant cloths and shaved 
his head, in the memory of his father, and he liked to wear sandals
made of tire. When he was invited to Harvard University
in 1958, he criticized the visa requirement of swearing that he
would refrain from subversive actions.
Like other prominent mathematicians, he opposed to the French 
colonial war in Algeria, 1954-62.
In 1966 he was awarded the {\em Fields Medal,} which he was supposed to 
receive in Moscow, but he did not show up, referring to two Russian 
writers who were arrested.\footnote{IH\'{E}S director L\'{e}on 
Motchane received the Fields Medal at the {\it International Congress 
of Mathematicians} in Moscow, on Grothendieck's behalf.}
Still, he did visit Eastern European countries
at other occasions, and his ideas had a remarkable influence among
Russian mathematicians, like Vladimir Drinfeld, Maxim Kontsevich, 
Yuri Manin and Vladimir Voevodsky. 

Meanwhile the students movement of the 1960s gained more and 
more momentum, and culminated in May 1968 in Paris.
Grothendieck was strongly impressed, but he
found himself on the wrong side --- he felt attracted to 
the r\^{o}le of an outlaw, not to the establishment.
He sympathized with the movement, which involved in part
anarchist ideas, but he did not attend public rallies; 
also in this respect he followed his own path.

In 1967 Grothendieck received a request from Hanoi, 
asking for literature about Algebra and Algebraic Geometry. 
He had not been aware that there was mathematical research going on
in North Vietnam, even during the worst period of the second
Vietnam War, and he provided as much material as he could. Moreover, 
Grothendieck decided to travel to North Vietnam himself to give
lectures. After a first part in Hanoi, the cluster bombing 
by the US Air Force intensified so much\footnote{A week
after Grothendieck's arrival, the campus of the Hanoi 
Polytechnical Institute was hit 
by delay-action bombs, which killed two mathematicians.} 
that he and his Vietnamese attendees 
(among them Ta Quang Buu, mathematician and Minister of 
Higher Education and Technology, who frequently asked questions)
moved to a hidden place in the forest to continue the
lecture (during the breaks, Grothendieck went to a near-by
river to wash his cloths). In 2013 Neal Koblitz, a mathematician
from the USA, visited this place, and he was intrigued by the
fact that the course given there could have been presented as
well at Harvard University, where Koblitz had been studying at that
time \cite{Vietnam1}. After his return to France, Grothendieck gave 
talks about his visit and wrote a detailed report \cite{Vietnam2}, 
which informed the world about the mathematical community in North 
Vietnam. While he described the state as somewhat over-regulated, 
his report firmly expresses his sympathy for the underdogs 
of this destructive war, which lasted for 30 years in total,
and left nearly 4 millions of people killed.
\begin{figure}
\begin{center}
\includegraphics[angle=0,width=.45\linewidth]{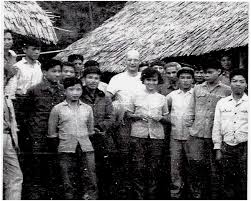}
\hspace*{1cm}
\includegraphics[angle=0,width=.45\linewidth]{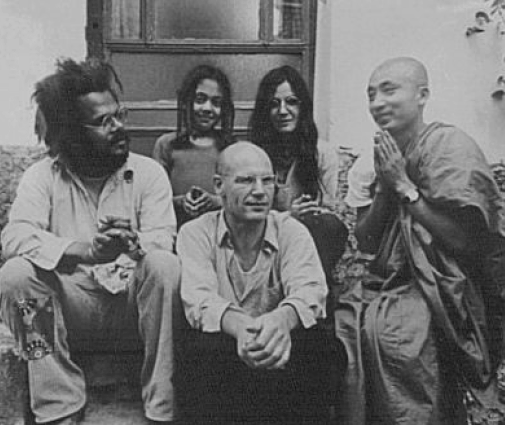}
\end{center}
\vspace*{-2mm}
\caption*{\it On the left: Grothendieck in the Vietnamese rain forest, 
100 km from Hanoi. On the right: a picture of 1975, characteristic for
his new life style.} 
\vspace*{-2mm}
\end{figure}

When Grothendieck quit IH\'{E}S in 1970, the reason he gave was
that he had discovered that his institute received funds for the
French military. Actually this was known before, and many
scientists objected, so in 1969 an agreement was reached
to stop this practice. However, only one year later the agreement
was broken. Grothendieck tried to convince his colleagues to resign
in protest, but it was only him who really did so.
Although this was about a minor fraction of the IH\'{E}S budget,
it was an ethical problem for him --- for instance,
his Ph.D.\ advisor Schwartz had been working hard to transform
the \'{E}cole Polytechnique from a military to a civilian
orientation. Research with military purposes was not
acceptable for Grothendieck, and this was also a reason why physics
was suspicious to him, keeping in mind Hiroshima and Nagasaki.\footnote{On 
the other hand, although his mathematical style was very abstract, some 
of his concepts did propagate into theoretical physics, in particular
in constructive field theory, which was elaborated in the 1970s. A
later example is {\em non-commutative geometry,} which was --- on 
the formal side --- strongly developed by Alain Connes,
who adopted ideas by Grothendieck. It became a wide-spread
fashion in theoretical physics in the late 1990s.
Yet another example is the {\em Atiyah-Singer Index Theorem,} which 
relates the zero modes of a chiral Dirac operator to the topological
charge of its gauge background. Its derivation involved Grothendieck 
Groups, which emerged from his new proof of the Riemann-Roch Theorem,
and which became later a point of departure for K-Theory.\label{NC}}

However, there might have been other reasons involved in this 
radical turning point of his life: conflicts with colleagues at 
IH\'{E}S, a decrease in creativity, and the consciousness that
his ambitious goals would never be completed.\footnote{He had 
outlined a monumental program to write a series entitled {\it 
\'{E}l\'{e}ments de G\'{e}om\'{e}try Alg\'{e}brique} in 13 volumes.
``Only'' four volumes appeared, comprising about 1800 pages.}
Was there also some burnout or mid-life crises involved?

In any case, Grothendieck changed his life style, he got
separated from his wife Mireille, and opened {\it communes},
first in Paris and later in Southern France. There he lived
with a variety of people, at times three of his children
were among them, and political meetings were held. Meanwhile
he lectured on a temporary basis, first at the 
{\it University Paris-Sud} in Orsay, and then in the 
{\it College de France.} In his courses he took the 
opportunity to discuss also issues like the
threat by nuclear weapons. This attracted a broad audience,
but the {\it College} direction was not amused, and denied him a
permanent position,  even though he was one of the most famous
mathematicians in the world.

In 1973 he moved back to the University of Montpellier (although its 
Mathematics Department did still not match his standard)
and gave lectures on all levels. He was friendly to his students, 
who dubbed him {\it Alexandre le Grand,}
he distributed organic apples, and gave inspiring
courses. He did not run a highly ambitious seminar anymore, but
he still had several Ph.D.\ students (and he got angry when the 
{\it Springer} Publishing House declined publishing a thesis).
Still he led excellent research, but the French research agency
CNRS only provided marginal support.

From 1973-79 he lived in the tiny village Olmet-et-Villecun, 
50 km from Montpellier, in a simple house without electricity
(he used kerosene lamps to work at night).
He did not hesitate to give shelter to homeless people.
Generally, his home was open for everybody, and it became a 
meeting place for all kind of people, including the {\em hippie 
movement.} In 1977 it was raided by the police,
which were looking for anything possibly illegal.
All they found was a Japanese citizen, who 
was staying there, and whose French visa had expired. 
He was an peaceful person, who had studied mathematics,
but at that time he was a Buddhist monk. Half a year later 
(when the monk had long left France), Grothendieck was actually accused 
for giving shelter and food to a foreigner ``in an irregular situation''. 
He defended himself with a passionated speech, and many mathematicians 
gave him public support, but he got convicted to a heavy
fine and a six-months suspended sentence.

\section{Environmental and peace movement}

Meanwhile Grothendieck questioned intensively the sense of
scientific research --- he reported that in many discussions,
nobody could really give a reason for it.
He got more and more concerned about ecological
problems and militarism, in particular the danger 
of a nuclear war. He was 
convinced that everyone, who was given the relevant information,
would follow his logical arguments and agree with his conclusions, 
and that he had a mission to spread this message. 

When he was invited to lecture at a Summer School in Montreal, 
he accepted under the condition that he could not only lecture 
about mathematics, but also about the threats to humanity.
In fact, some young mathematicians followed
his ideas, and became activists as well. He also gave double
lectures in the USA, where he further supported the rights of
Native Americans. A Ph.D.\ student named Justine Skalba was 
particularly excited about his charisma and intelligence. She
followed him to France as his partner for a few years
(their son was born after the relation had ended, but he 
later did a Ph.D.\ in mathematics). Justine
remembers a rally in Avignon that they attended together; 
when it was harassed by the police, Grothendieck
knocked down two policemen, and got arrested.

Together with another two prominent French mathematicians, 
Claude Chevalley and Pierre Samuel, he founded
a group called {\it International Movement for the Survival
of the Human Race.} It published the magazine {\it Vivre,} or 
later {\it Survivre et Vivre,} with emphatic calls for peace and
against pollution, discussions of the impact of science, and a 
critical view of the consumerism-minded society (in part inspired 
by the philosopher Herbert Marcuse). He wrote a considerable number 
of articles for this magazine, which appeared from 1970-75; copies 
are available in the internet \cite{circle}.\footnote{It goes 
without saying that this movement was also confronted with dismissive
reactions; {\it e.g.}\ Ref.\ \cite{Cartier} describes it
as a ``dooms-day sect'', which was ``obsessed by pollution''.}

When Grothendieck attended the {\it International Congress of Mathematicians}
in Nice, 1970, he installed a desk to distribute this magazine, together
with his eldest son Serge (from a premarital relation), and tried 
to recruit new members to his movement (with limited success). 
Dieudonn\'{e}, who  was responsible for the event, stubbornly 
objected, until they moved the desk outside the building, but there 
they got in trouble with the police. 

At a Summer School 1972 in Antwerp, Serre gave the opening speech. 
Grothendieck vociferously interrupted his former IH\'{E}S colleague, 
to speak out against NATO, which had sponsored this event.\footnote{As
a reaction, a NATO representative, who had intended to join the
Summer School for a public debate, backed off. Subsequently
Grothendieck was blamed for having done some kind of ``damage''.} 
He did not hesitate to be provocative {\it (l'enfant terrible)}, 
even if this led to resentment with long-term friends and 
collaborators.

Of course he was confronted with the reproach of overdoing
it in an immature manner. However, for instance
a talk of two hours, including an extensive discussion with the
audience, that he gave at CERN in 1972 (now accessible
on YouTube) sounds calm and thoughtful: he was aware, of course,
that CERN does not focus on nuclear research (in contrast
to its name). He explained why he took distance 
from the scientific community, with its competition and
pressure to publish, which are unjust and unfavorable for
creativity, and which keeps researchers working without ever
wondering for what reason. He also recalled mathematicians who
had committed suicide. He further pointed out why he now 
considered actions against the threats to humanity --- like 
nuclear weapons --- far more important.

\begin{figure}
\begin{center}
\includegraphics[angle=0,width=.25\linewidth]{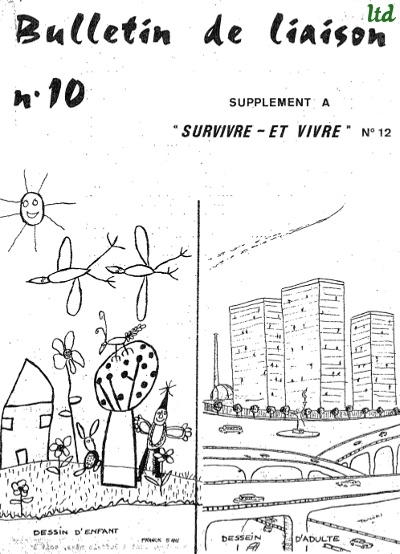}
\hspace*{2cm}
\includegraphics[angle=0,width=.25\linewidth]{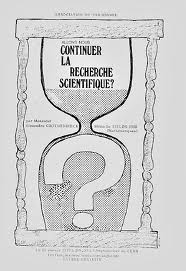}
\end{center}
\vspace*{-2mm}
\caption*{\it On the left: a page of the magazine ``Survivre et Vivre'':
it compares the drawing by a child and by an adult.
On the right: announcement of Grothendieck's talk at CERN, 1972, 
entitled ``Are we going to continue the scientific research?''.}
\end{figure}

Despite his conviction and arguments, the use of his
academic reputation and his rhetorical skills, his group remained
small. It mostly attracted people who had similar views before,
and around 1973 a trend of dissolution set in. Grothendieck
was disappointed and considered his efforts as a failure.
He concluded that people, even scientists, were blind to the
dangers to the world, and do not behave rationally.

Other members of this movement, like Samuel (who was an editor
of  {\it Survivre et Vivre} until 1973), patiently carried on
efforts along these lines.
From today's perspective, these actions appear as pioneering 
work for the peace and environmental movement, which later 
became influential --- to some extent --- in Europe and beyond.
Also within science, ecological concerns were later acknowledged,
{\it e.g.}\ with the 1995 Chemistry Nobel Prize for demonstrating
the danger to the Earth's ozone layer. Today, for instance global 
warming due to human activities is only disbelieved by some people 
who are at odds with science. So we could ask if Grothendieck's
activism in the early 1970s was really immature, or if he was
rather ahead of his time?

\section{Spiritualism and isolation}

Frustrated about his modest success as an activist, Grothendieck 
slowed down his appeals to the public. He kept writing long 
manuscripts, with magnitudes of 1000 pages, like {\it La Longue 
Marche \`{a} travers la th\'{e}orie de Galois,}\footnote{Grothendieck
had a special esteem for \'{E}variste Galois, whom he called his
{\it fr\`{e}re de temp\'{e}rament} \cite{Recoltes}. In fact the 
two mathematical geniuses had a number of points in common: 
a push for a new level of abstraction and a major interest
in the {\em relation} among mathematical objects; now 
the term ``Grothendieck's Galois Theory'' is used. Moreover, both 
had an early end of their career (Galois' case was far more extreme), 
and they were radical activists for what each one considered, 
in his epoch, as an urgent progress of society.}
{\it A la Pursuite de Champs} and {\it Esquisse d'un Programme}
with ideas for future mathematics. Indeed, that {\it Programme} 
was worked out to a large extent by the young mathematicians
Leila Schneps and Pierre Lochak, who were impressed by its
farsighted vision. They contacted Grothendieck, who suddenly
expressed his interest in {\em physics} and asked for literature
about it \cite{Jackson} (although he regretted its lack of rigor). 
Later they also initialized the {\it Grothendieck Circle,} which 
created an informative web page \cite{circle}, and Schneps edited an 
overview over Grothendieck's mathematical achievements \cite{Schneps}.

In the period 1983-88, Grothendieck wrote a stylistically brilliant 
book entitled  {\it R\'{e}coltes et Semailles} \cite{Recoltes}, where he 
reviews his life and work, supplemented by all kinds of elements,
like love poems (in German) and (sometimes critical) comments 
on the mathematical community and former colleagues.
In Section 2.20 he addressed modern physics. 
From a mathematical perspective, he did not consider
Einstein's Theory of Relativity very interesting, although he
appreciated its importance for our paradigm of space-time.
Mathematically, however, he described the transition from Newton's 
Theory to Relativity like a change from one French dialect to another, 
whereas Quantum Theory is like a transition to Chinese. 
This he did find interesting, regarding his deep understanding 
of a {\em point,} and he mentions an intuitive similarity to his 
concept of {\em toposes.} We add a comment in Appendix A.

In 1988 he was supposed to receive the prestigious {\em Crafoord Prize} 
by the Royal Swedish Academy of Sciences, together with his former
student Deligne, but Grothendieck declined.
In a polite letter \cite{Crafoord} he explained his reasons: 
first he did not need money, and about the importance of his 
work, time and offspring would decide, not honors.
He adds that such prizes are constantly given to the
wrong people, who do not need further wealth
nor glorification. He asks whether this ``superabundance 
for some'' is not provided ``at the cost of the needs of others''?
Finally, he points out that agreeing to ``participate in the game 
of {\it prizes}'' would imply his ``approval to the spirit
$\dots$ of the scientific world'', where 
ethics has ``declined to the point that outright theft among
colleagues (especially at the expenses of those who are in no 
position to defend themselves) has nearly become a general rule''.

Also in 1988 he retired from Montpellier University, and in 1991 even 
from {\em society;} he broke off contacts with almost everybody,
including his family. He withdrew to a modest life in a hamlet
in the French Pyrenees, not far from Vernet Camp 
(the redoubtable camp, where his father had been imprisoned 
before being deported to Auschwitz). He still wrote the mathematical 
program {\it Les D\'{e}rivateurs} (about 2000 pages), which he 
handed over to a friend.
On the other hand, he once burned a huge amount of notes, 
letters and other documents, one estimates 25\,000 pages. His 
main interest now shifted to spiritualism and meditation,
and he entered the final, Steppenwolf-like phase of his life. 

For quite some time, since the 1970s, he was strongly interested in 
Buddhism. There are hints that this helped him to relax from the 
pressure of productivity, and to improve the relation with his 
ex-wife Mireille. He was a strict vegetarian and received Buddhist 
teachers. He was also fascinated by the symbols of Yin and Yang
\cite{Recoltes}, and characterized his style of research as 
Yin.\footnote{Deligne describes a proof by Grothendieck as a lengthy
sequence of trivial steps, ``nothing seems to happen, but yet
at the end a highly non-trivial theorem is there''. This is in
contrast to Serre's Yang-style of striving for a solution in one 
strike \cite{RisingSea}.}
Later, however, he moved on to a mystic
and unconventional form of Christianity.
He spent a period of starving, which endangered his health.
He got much interested in dreams, which he considered
the messengers of spiritual wisdom, and he studied
Freud's interpretation. 

As his main activity, he kept on writing; daily he spent many 
hours typing about his mystic experiences and 
ideas, which led again to thousands of pages. Although he assumed 
his visions to be relevant for the future society, he did not
want to publish these notes. In 2010 a bizarre (but well formulated)
hand-written letter appeared, where he even requested the removal of 
all his works from the libraries.

\begin{figure}
\begin{center}
\includegraphics[angle=0,width=.25\linewidth]{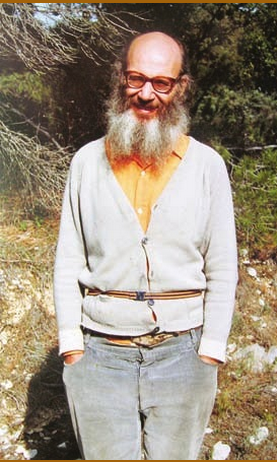}
\end{center}
\vspace*{-2mm}
\caption*{\it Grothendieck in his older days.} 
\vspace*{-2mm}
\end{figure}
Only very few people knew where he lived and promised not to spread
this information. He did not have a postal address nor a telephone,
let alone internet, and he did not receive uninvited visitors.
At last the world learned that he has died, on November 13, 2014, 
in the village of Saint-Lizier in the French Pyrenees --- 
Alexander Schurik Grothendieck, rest in peace.

\appendix

\section{Points in physics}

When Grothendieck writes that he finds Einstein's Relativity 
mathematically ``banal'' \cite{Recoltes}, he includes General 
Relativity. That is based on mathematics of the 19th century,
in particular Differential Geometry, which is actually
non-trivial (by common perception). 
Nevertheless, his statement becomes plausible if it addresses
a specific, fundamental understanding of geometry, in particular
the very nature of a {\em point;} we have mentioned before 
that this issue has haunted Grothendieck since his youth. In this
regard, he considers Quantum Mechanics far more interesting.

In Ref.\ \cite{Recoltes}, p.\ 69, he writes: 
``And these {\it probability clouds,} which replace the certain
material particles that we had before, remind me strangely of
the elusive {\it open neighborhoods} which inhabit the toposes,
such evanescent ghosts, which surround fictitious 
{\it points,} which keep on attaching themselves, 
in contrast to a recalcitrant imagination''.
 
The state of a quantum mechanical particle is given by a 
(time-dependent) vector in a Hilbert space, in Paul Dirac's 
notation $| \psi (t) \rangle$, and the position
eigenstates $| x \rangle$ form a basis. The scalar product
$\psi (t,x) = \langle x | \psi (t) \rangle$ is the particle's 
wave function, and $| \psi (t,x)|^{2}$ its ``probability cloud''.

Hence standard Quantum Mechanics is still formulated
in a standard coordinate space with a continuum of sharp points, 
where the particle wave functions are accommodated. The spatial 
resolution is not limited in principle, if sufficiently large 
momenta are available to resolve it. So a possible analogy
to {\em toposes} could rather refer to {\em phase space,} 
where points do have a conceptually limited resolution --- given by 
Heisenberg's Uncertainty Relation --- hence there are only 
{\em fuzzy} points.

On the other hand, angular momenta only take sharp, discrete values.
A deep understanding of this interplay between discrete and continuous, 
sharp and fuzzy, could have been something for Grothendieck's taste.

There have been attempts to employ {\em toposes} in an 
unconventional formulation of Quantum Mechanics, with the hope 
to alleviate problems with the interpretation of measurements on 
quantum systems; for reviews see Refs.\ \cite{toposQM}.

We still do not have a convincing merger of Quantum Physics
with General Relativity, but a Lorentz covariant formulation
(hence a conciliation with Special Relativity) is accomplished
by Quantum Field Theory. But again, the usual formulation,
which encompasses in particular the (tremendously successful)
Standard Model of Particle Physics, employs a simple Minkowski 
space. Thus the fields are still functions of sharp points
in space-time --- this does still not seem exciting in view of 
Grothendieck's particular motivation, at least at first sight.

One might object that some regularization of high momentum contributions 
is required, to suppress the omnipresent {\em ultraviolet singularities.}
This corresponds to some kind of truncation of short distances, 
{\it i.e.}\ a somehow granular space-time (its structure depends on
the regularization scheme).\footnote{The scheme that works at finite
interaction ({\it i.e.}\ beyond perturbation theory) performs indeed
a reduction to a ``lattice'' of discrete (but sharp) space-time points.}
In theories like Quantum Chromodynamics
(the sector of the Standard Model that describes the strong interaction),
this truncation can be fully removed at the end of the calculation,
{\it i.e.}\ one extrapolates to the {\em continuum limit,}
hence it is just a mathematical trick. Also the electroweak
sector was shown to be {\em renormalizable,} which allows to
remove the cutoff at the end. On the other hand, in
the Higgs sector of the Standard Model,
a complete removal of the truncation would also remove
the interactions, so the Higgs field becomes free and does no longer
do its job of providing the elementary particle masses. Usually physicists 
don't worry much about this property (which is known as {\it triviality}),
since a huge momentum cutoff, extremely far above the experimentally 
accessible regime, is sufficient to justify the observed mass of the 
Higgs boson.\footnote{It is now getting popular to interpret 
the Standard Model as a low energy effective theory, valid up to some 
energy range above the scale of experiments, which is sufficient for 
practical purposes.}
However, from a fundamental geometric point of
view, one might pay more attention to this aspect.

Based on reports of people who talked to Grothendieck in his old
days, he seemed interested in the question if the constants of
Nature are related by {\em rational} ratios. Usually we do not assume that
({\it e.g.}\ we do not have reasons to expect the fine structure
constant $\alpha = e^{2}/(4 \pi \epsilon_{0} \hbar c) \simeq
1/(137.036)$, or the ratios of particle masses, to be rational), but
it is certainly the case for the electric particle charges. 
Dirac gave an explanation for this property, but it
requires the existence of at least one magnetic monopole, for
which we do not have any evidence. In the Standard Model and
some (though not all) of its variants --- in particular incorporating
neutrino masses --- that property can also be deduced from the 
theoretical requirement of gauge anomaly cancellations.

In the framework of the diverse attempts to achieve 
compatibility of Quantum Theory also with General Relativity, 
{\it i.e.}\ to include also gravity, quite general 
arguments suggest an extension to a {\em pure space-time 
uncertainty relation,} which should be manifest
at extremely short distances, of the order of the Planck scale
($\approx 10^{-35}$ m), see {\it e.g.}\ Ref.\ \cite{DFR}.
That corresponds to a {\em non-commutative geometry,} where the 
coordinates in independent directions are given by Hermitian 
operators, which do not commute. We speculate that this could
have attracted Grothendieck's interest.

We have mentioned in footnote \ref{NC} that the formal
aspects of physics in a non-commutative space have been elaborated
rigorously by Connes. He did apply some concepts by Grothendieck, 
but regarding {\it toposes} his comprehensive book on this subject,
Ref.\ \cite{Connes}, only contains the remark: ``One could base 
this extension of topology on the notion of toposes due to Grothendieck. 
Our aim, however, is to establish contact with the powerful tools 
of functional analysis such as positivity and Hilbert space 
techniques, and with K-theory.''

From the physical perspective we add that Quantum
Field Theory --- formulated in a non-commutative space --- 
is plagued by severe obstacles:
first we cannot instal fields for the non-Abelian 
gauge groups $SU(2)$ and $SU(3)$, which belong to the
Standard Model. Next a non-commutative space-time entails 
{\em non-local interactions} --- which also occur in String Theory 
--- and which raise questions regarding the principle of causality. 
That might be acceptable if these non-local effects were restricted 
to tiny ranges (like the Planck scale), but for interacting quantum
fields, a non-commutative space-time further gives rise to a new type 
of singularity in the {\em infrared regime.} Hence quantum effects 
are expected also at very {\em long} distances, even if one modifies 
the geometry only within a tiny range. This theoretical phenomenon 
--- known as ``ultraviolet-infrared mixing'' --- has prevented 
a valid confrontation with particle phenomenology, which could confirm
or constrain a space-time non-commutativity in Nature, and thus the
existence of fuzzy points, which bear a resemblance to the
open neighborhoods in the Grothendieck toposes. \\

{\it A shorter version of this article has been published
in Spanish \cite{CIENCIAS}.}

\end{document}